\input amstex
\documentstyle{amsppt}
\magnification1200
\NoBlackBoxes
\pagewidth{6.5 true in}
\pageheight{9.25 true in}
\document
\topmatter 
\title  Quantum Unique Ergodicity for $SL_2({\Bbb Z})\backslash {\Bbb H}$ 
\endtitle
\author
 K. Soundararajan
\endauthor 
\address
{Department of Mathematics, 450 Serra Mall, Bldg. 380, Stanford University, 
Stanford, CA 94305-2125, USA}
\endaddress
\email 
ksound{\@}stanford.edu
\endemail
\thanks   The author is partially supported by the National Science Foundation (DMS 0500711).
\endthanks
\endtopmatter

 \document
\head 1.  Introduction \endhead 

\def\lam{\lambda}
\noindent An interesting problem in number theory and 
quantum chaos is to understand the distribution of 
Maass cusp forms of large Laplace eigenvalue for the modular surface
$X= SL_2({\Bbb Z})\backslash {\Bbb H}$.  Let $\phi$ denote 
a Maass form of eigenvalue $\lambda$, normalized so that 
its Petersson norm $\int_{X} |\phi(z)|^2 \frac{dx \ dy}{y^2} $ 
equals $1$.  Zelditch [14] has shown that as $\lam \to \infty$, 
for a {\sl typical} Maass form $\phi$ the associated 
probability measure $\mu_{\phi} := 
|\phi(z)|^2 \frac{dx \ dy}{y^2}$ tends to the uniform 
distribution measure $\frac 3{\pi} \frac{dx\ dy}{y^2}$.  
This result is known as ``Quantum Ergodicity."  The widely studied 
Quantum Unique Ergodicity conjecture of Rudnick and 
Sarnak [10]  asserts that as $\lam \to \infty$, for {\sl every} Maass 
form $\phi$ the measure $\mu_{\phi}$ approaches the  
uniform distribution measure.  In studying this conjecture,  
it is also natural to restrict to Maass forms that are eigenvalues of 
all the Hecke operators; it is expected that the spectrum of the 
Laplacian on $X$ is simple so that this condition would 
automatically hold, but this is far from being proved. Using methods from ergodic 
theory, Lindenstrauss [6] has made great 
progress towards the QUE conjecture for such Hecke-Maass forms.  Namely, 
he has shown that the only possible weak-$*$ limits of the measures $\mu_{\phi}$ 
are of the form $c \frac{3}{\pi} \frac{dx \ dy}{y^2}$ where $c$ is 
some constant in $[0,1]$.  In other words, Lindenstrauss establishes 
QUE for $X$ except for the possibility that for some infinite subsequence 
of Hecke-Maass forms $\phi$ some of the $L^2$ mass of $\phi$ could 
``escape'' into the cusp of $X$.  In this paper we eliminate 
the possibility of escape of mass, and together with Lindenstrauss's 
work this completes the proof of QUE for $X$.  
 
 The results of Zelditch and Lindenstrauss are in fact  
 stronger than we have indicated above.   Given a Maass form 
 $\phi$ on $X$, Zelditch defines the ``micro-local" lift ${\tilde \phi}$ of $\phi$ to 
 $Y= SL_2({\Bbb Z}) \backslash SL_2({\Bbb R})$.  This micro-local 
 lift defines a measure on $SL_2({\Bbb Z})\backslash SL_2({\Bbb R})$ 
 with two important properties:  First, as the eigenvalue tends to infinity, the projection of 
 the measure from $Y$ to the surface $X$ approximates the measure $\mu_{\phi}$ 
 given above.  Second, as the eigenvalue tends to infinity any weak-$*$ limit of 
 these measures on $Y$ is invariant under the geodesic flow on $Y$.  
 Zelditch's result then asserts that for a full density subsequence of eigenfunctions, 
 the associated micro-local lifts get equidistributed on $Y$.  Lindenstrauss's 
 result is that any weak-$*$ limit of the lifts arising from Hecke-Maass forms
  is a constant $c$ (between $0$ and $1$) 
 times the normalized volume measure on $Y$.  We remark that 
 the analog of quantum unique ergodicity with Eisenstein 
 series in place of cusp forms has been treated by Luo and 
 Sarnak [8] in the modular surface version, and by Jakobson [5]
 for the corresponding micro-local lifts.   For more complete 
 accounts of the quantum unique ergodicity 
 problem the reader may consult [7], [9], [11], [13] and references therein; 
 a comprehensive introduction 
 to the theory of Maass forms is provided in [4].

\proclaim{Theorem 1}  Let $\phi$ be a Hecke-Maass cusp form 
for the full modular group $SL_2({\Bbb Z})$, normalized to 
have Petersson norm $1$.  Let $\tilde \phi$ denote the 
micro-local lift of $\phi$ to $Y= SL_2({\Bbb Z}) \backslash SL_2({\Bbb R})$, and 
let $\tilde{\mu_{\phi}}$ denote the corresponding measure on $Y$.   The 
normalized volume measure on $Y$ is the unique weak-$*$ limit 
of the measures $\tilde{\mu_{\phi}}$.  In particular, 
for any compact  subset $C$ of a fundamental domain for $SL_2({\Bbb Z})\backslash {\Bbb H}$ 
we have, as $\lam \to \infty$, 
$$ 
\int_{C} |\phi(x+iy)|^2 \frac{dx \ dy}{y^2} = \int_{C} \frac{3}{\pi} \frac{dx \ dy}{y^2} + o(1). 
$$ 
\endproclaim 

Theorem 1 is a consequence of the following result 
which estimates how much of the mass of $\phi$ can 
be present high in the cusp.

\proclaim{Proposition 2}   Let $\phi$ denote a Hecke-Maass 
cusp form for $SL_2({\Bbb Z})$ with eigenvalue $\lambda$,  
and normalized to have Petersson norm $1$.  For $T\ge 1$, we have 
$$ 
\int\Sb |x|\le \frac 12 \\ y \ge T \endSb |\phi(x+iy)|^2 \frac{dx \ dy}{y^2} 
\ll \frac{\log (eT)}{\sqrt{T}}.
$$
\endproclaim 

We remark that by entirely different methods Holowinsky and Soundararajan 
([2], [3], [12]) have settled the holomorphic analog of QUE for $X$; it is not clear 
how to adapt Lindenstrauss's methods to that setting.  Their 
methods have the advantage of yielding explicit estimates for the 
rate of convergence to uniform distribution; it is not clear 
how to obtain such a rate of convergence in Theorem 1.  However 
the works of Holowinsky and Soundararajan  use in an essential 
way Deligne's bounds for the Hecke eigenvalues of holomorphic 
modular forms; the analog of these bounds for Maass forms 
remains an important open problem.  

While we have restricted ourselves 
to the full modular group, our argument would apply also to all congruence 
subgroups.  Thus QUE for Maass forms, and its holomorphic analog, are 
now known for non-compact arithmetic quotients of ${\Bbb H}$.  In the 
case of compact arithmetic quotients, Lindenstrauss's work establishes
QUE for Maass forms; the analog for holomorphic forms 
remains open.  

Our proofs of Theorem 1 and Proposition 2 exploit the particular multiplicative 
structure of the Hecke-operators.  We say that a function $f$ is {\sl Hecke-multiplicative}
if it satisfies the Hecke relation
$$ 
f(m)f(n) = \sum_{d|(m,n)} f(mn/d^2), 
$$ 
and $f(1)=1$.  
The key to establishing Proposition 2 is the following result on 
Hecke-multiplicative functions. 

\proclaim{Theorem 3}  Let  $f$ be a Hecke-multiplicative 
function.  Then for all $1\le y\le x$ we have
$$ 
\sum_{n\le x/y} |f(n)|^2 \le 10^8 \Big(\frac{1+\log y}{\sqrt{y}} \Big) \sum_{n\le x} |f(n)|^2.
$$ 
\endproclaim 

It is noteworthy that Theorem 3 makes no assumptions on the size of the 
function $f$.   Hecke-multiplicative functions satisfy $f(p^2) = f(p)^2-1$, 
so that at least one of $|f(p)|$ or $|f(p^2)|$ must be bounded away from zero; 
this observation plays a crucial role in our proof.  
We also remark that apart from the $\log y$ factor, Theorem 3 
is best possible:  Consider the Hecke-multiplicative function $f$ defined by $f(p)=0$ for 
all primes $p$.  The Hecke relation then mandates that $f(p^{2k+1}) =0$ and 
$f(p^{2k}) = (-1)^k$.  Therefore, in this example, $\sum_{n\le x} |f(n)|^2 = \sqrt{x} + O(1)$ 
and $\sum_{n\le x/y} |f(n)|^2 = \sqrt{x/y} + O(1)$.  

The argument of Theorem 3 can be generalized in several ways.  For example 
one could obtain an analogous result with $|f(n)|$ in place of $|f(n)|^2$.  
Moreover one could consider multiplicative functions $f$ arising 
from Euler products of degree $d$.  By this we mean that for each prime $p$ 
there exist complex numbers $\alpha_j$ ($j=1,\ldots, d$) with $|\alpha_1 \cdots \alpha_d|=1$ 
and $\sum_{k=0}^{\infty} f(p^k) x^k = \prod_{j=1}^{d} (1-\alpha_j x)^{-1}$; the 
case $d=2$ corresponds to our Hecke-multiplicative functions above.   For these functions, 
one of $|f(p)|$, $\ldots$, $|f(p^d)|$ must be bounded away from zero, and 
exploiting this we may establish an analog of Theorem 3.   


{\bf Acknowledgments.}  I am grateful to Peter Sarnak for encouragement and 
some helpful suggestions, and to Roman Holowinsky for a careful reading.

\head 2.  Deducing Theorem 1 and Proposition 2 from Theorem 3 \endhead 
\def\lam{\lambda}

\demo{Proof of Proposition 2}  Let $\phi$ be a Maass form of eigenvalue $\frac 14+ r^2$ for 
the full modular group, normalized to have Petersson norm $1$.   
We suppose that $\phi$ is an  eigenfunction of all the Hecke operators, 
and let $\lam(n)$ denote the $n$-th Hecke eigenvalue.  
Recall that $\phi$ has a Fourier expansion of 
the form
$$ 
\phi(z) =C \sqrt{y} \sum_{n=1}^{\infty} \lam(n) K_{ir}(2\pi ny) \cos (2\pi nx),
$$ 
or 
$$ 
\phi(z) = C \sqrt{y} \sum_{n=1}^{\infty} \lam(n) K_{ir}(2\pi ny) \sin(2\pi nx), 
$$ 
where $C$ is a constant (normalizing the $L^2$ norm), $K$ denotes the usual $K$-Bessel 
function, and we have $\cos$ or $\sin$ depending 
on whether the form is even or odd.  

Using Parseval we find that 
$$ 
\int\Sb |x|\le \frac 12 \\ y \ge T \endSb |\phi(x+iy)|^2 \frac{dx \ dy}{y^2} 
= \frac{C^2}{2} \int_{T}^{\infty} \sum_{n=1}^{\infty} |\lam(n)|^2 |K_{ir}(2\pi ny)|^2 \frac{dy}{y}.
$$
By a change of variables we may write this as 
$$ 
\frac{C^2}{2} \sum_{n=1}^{\infty} |\lam(n)|^2  \int_{nT}^{\infty} |K_{ir}(2\pi t)|^2 \frac{dt}{t} =
\frac{C^2}{2} \int_{1}^{\infty} |K_{ir}(2\pi t)|^2 \sum_{n\le t/T} |\lam(n)|^2 \frac{dt}{t}.
$$ 
Appealing to   Theorem 3 this is 
$$ 
\align
&\ll \frac{\log eT}{\sqrt{T}} \frac{C^2}{2 }  \int_1^{\infty} |K_{ir} (2\pi t)|^2 \sum_{n\le t} |\lam(n)|^2 \frac{dt}{t} \\
&=\frac{\log eT}{\sqrt{T}} \int\Sb |x|\le \frac 12 \\ y \ge 1 \endSb |\phi(x+iy)|^2 \frac{dx \ dy}{y^2} 
\ll \frac{\log eT}{\sqrt{T}}, \\
\endalign
$$ 
since the region $|x|\le \frac 12$, $y\ge 1$ is contained inside a 
fundamental domain for $SL_2({\Bbb Z})\backslash {\Bbb H}$.  This proves Proposition 2. 
\enddemo 

\demo{Proof of Theorem 1} As remarked in the introduction, Lindenstrauss has 
shown that any weak-$*$ limit of the micro-local lifts of Hecke-Maass forms 
is a constant $c$ (in $[0,1]$) times the normalized volume measure on $Y$.  
Projecting these measures down to the modular surface, we see that 
any weak-$*$ limit of the measures $\mu_\phi$ associated 
to Hecke-Maass forms is of the shape $c \frac{3}{\pi } \frac{dx \ dy}{y^2}$.
Theorem 1 claims that in fact $c=1$, and there is no escape of 
mass.  If on the contrary $c<1$ for some weak-$*$ limit, then we have a sequence of 
Hecke-Maass forms $\phi_j$ with eigenvalues $\lam_j$ tending to infinity 
such that for any fixed $T\ge 1$ and as $j\to \infty$
$$ 
\int\Sb z \in {\Cal F}\\ y\le T \endSb |\phi_j(z)|^2 \frac{dx \ dy}{y^2} = (c+o(1)) \frac{3}{\pi} \int\Sb z\in {\Cal F} \\ y\le T\endSb \frac{dx \ dy}{y^2}  
= (c+o(1)) \Big( 1- \frac{3}{\pi T} \Big); 
$$ 
here ${\Cal F} = \{ z= x+iy: |z| \ge 1,  \  -1/2\le x\le 1/2, y>0\}$ denotes the 
usual fundamental domain for $SL_2({\Bbb Z})\backslash {\Bbb H}$.  It follows 
that as $j\to \infty$ 
$$ 
\int\Sb |x| \le \frac 12 \\ y \ge T \endSb |\phi_j(z)|^2 \frac{dx \ dy}{y^2} 
=  1- c + \frac{3}{\pi T} c + o(1),
$$ 
but if $c<1$ this contradicts Proposition 2 for suitably large $T$.  

\enddemo 

\head 3.  Preliminaries for the proof of Theorem 3 \endhead 

\noindent Throughout the proof of Theorem 3 we shall work with a single 
value of  $x$.  Accordingly, we define for $1\le y\le x$ 
$$ 
{\Cal F}(y) = {\Cal F}(y;x) = \frac{\sum_{n\le x/y} |f(n)|^2}{\sum_{n\le x} |f(n)|^2}, 
$$ 
and our goal is to show that ${\Cal F}(y) \le 10^8 \log (ey)/\sqrt{y}$.
Note that ${\Cal F}(y)\le 1$ for all $y \ge 1$, and that ${\Cal F}$ is a decreasing function 
of $y$.  Thus in proving Theorem 3 we may assume that $y\ge 10^{16}$.  
We adopt throughout the convention that $f(t)=0$ when $t$ is not a natural
number, and that ${\Cal F}(y)=0$ if $y>x$.  

\proclaim{Lemma 3.1}  For any prime $p\le x$ we have 
$$ 
|f(p)| \le \frac{2}{{\Cal F}(p)^{\frac 12}}, 
$$ 
and, for $p\le \sqrt{x}$ 
 $$ 
 |f(p)| \le \frac{2}{{\Cal F}(p^2)^{\frac 14}} . 
 $$
\endproclaim 

\demo{Proof} Consider 
$$ 
|f(p)|^2 \sum_{n\le x/p} |f(n)|^2  = |f(p)|^2 {\Cal F}(p) \sum_{n\le x} |f(n)|^2. \tag{3.1} 
$$ 
Now $|f(p)f(n)| \le |f(pn)| + |f(n/p)|$ so that $|f(p)f(n)|^2 \le 2 (|f(np)|^2 + |f(n/p)|^2)$. 
Hence the LHS of (3.1) is 
$$ 
\le 2 \sum_{n\le x/p} (|f(np)|^2 + |f(n/p)|^2) \le 4 \sum_{n\le x} |f(n)|^2,
$$ 
and the first bound  follows. 

To see the second bound consider 
$$ 
|f(p^2)|^2 \sum_{n\le x/p^2} |f(n)|^2 = |f(p^2)|^2 {\Cal F}(p^2) \sum_{n\le x} |f(n)|^2.
\tag{3.2}
$$
The Hecke relations give $|f(p^2) f(n)| \le |f(p^2 n)| + |f(n)| + |f(n/p^2)|$ 
so that $|f(p^2) f(n)|^2 \le 3 ( |f(p^2 n)|^2 + |f(n)|^2 +|f(n/p^2)|^2)$.   Hence 
the LHS of (3.2) is $\le 9 \sum_{n\le x} |f(n)|^2$.  It follows that 
$$ 
|f(p^2)| \le \frac{3}{{\Cal F}(p^2)^{\frac{1}{2}}}, 
$$ 
and as $|f(p)|^2 \le |f(p^2)| + 1$ the result follows.  

\enddemo 

\proclaim{Proposition 3.2}  Let $d$ be a square-free number.  Then 
$$ 
\sum\Sb n\le x/y\\ d|n \endSb |f(n)|^2 \le \tau(d) \prod_{p|d} (1+|f(p)|^2) {\Cal F}(yd) \sum_{n\le x} |f(n)|^2, 
$$ 
where $\tau(d)$ denotes the number of divisors of $d$.  
Moreover
$$
\sum\Sb n\le x/y \\ d^2 |n \endSb |f(n)|^2 \le \tau_3(d) \prod_{p|d} (2+|f(p^2)|^2) {\Cal F}(yd^2) 
\sum_{n\le x} |f(n)|^2, 
$$ 
where $\tau_3$ denotes the $3$-divisor function (being the number of 
ways of writing $d$ as $abc$).
\endproclaim

\demo{Proof}   The Hecke relations give that $f(p)f(m) = f(pm)+f(m/p)$, and so $|f(pm)|\le |f(p)f(m)| +|f(m/p)|$. 
By induction  we find that 
$$ 
|f(md)| \le \sum_{ab=d} |f(a)| |f(m/b)|,
$$
so that 
$$ 
|f(md)|^2 \le \tau(d) \sum_{ab=d} |f(a)|^2 |f(m/b)|^2. 
$$ 
Summing over all $m\le x/(yd)$ we obtain that 
$$ 
\align
\sum\Sb n\le x/y\\ d|n \endSb |f(n)|^2 
&\le \tau(d) \sum_{ab=d} |f(a)|^2 \sum_{m\le x/(yd)} |f(m/b)|^2 
\\
&\le \tau(d) \sum_{ab=d} |f(a)|^2 {\Cal F}(yd) \sum_{n\le x} |f(n)|^2.\\
\endalign 
$$ 
The first statement follows. 

To prove the second assertion, note that the 
Hecke relations give that $f(m)f(p^2)$ equals 
$f(mp^2)$ if $p\nmid m$, $f(mp^2)+f(m)$ if $p\Vert m$, 
and $f(mp^2) + f(m) +f(m/p^2)$ if $p^2|m$.  In all 
cases we find that $|f(mp^2)| \le |f(p^2)f(m)| +|f(m)| +|f(m/p^2)|$.
By induction we may see that 
$$ 
|f(md^2)| \le \sum_{abc=d} |f(a^2)| |f(m/c^2)|.
$$ 
Therefore 
$$ 
|f(md^2)|^2 \le \tau_3(d) \sum_{abc=d} |f(a^2)|^2 |f(m/c^2)|^2.
$$ 
Summing over all $m\le x/(yd^2)$ we obtain that 
$$ 
\align
\sum\Sb n=md^2 \le x/y\endSb |f(n)|^2 &\le 
\tau_3(d) \sum_{abc=d} |f(a^2)|^2 \sum_{m\le x/(yd^2)} |f(m/c^2)|^2 
\\
&\le \tau_3(d) {\Cal F}(yd^2) \Big( \sum_{n\le x} |f(n)|^2\Big) 
\sum_{a|d} |f(a^2)|^2 \tau(d/a), 
\\
\endalign
$$
and the second statement follows.

\enddemo

Let ${\Cal P}= {\Cal P}(y)$ denote the set of primes in $[\sqrt{y}/2,\sqrt{y}]$.  The 
prime number theorem gives for large $y$ that $|{\Cal P}(y)| \sim \sqrt{y}/{\log y}$.  
In fact, using only a classical result of Chebyshev we find that for $y\ge 10^{16}$ we 
have $|{\Cal P}(y)| \ge \sqrt{y}/(2\log y)$ (see, for example, Dusart's thesis [1] which gives 
 more precise estimates).  The second 
bound of Lemma 3.1 gives that $|f(p)| \le 2/{\Cal F}(y)^{\frac 14}$.  Therefore,  we select 
$$ 
J = \Big[ \frac{1}{4\log 2} \log (1/{\Cal F}(y)) \Big] +3, 
$$ 
and partition ${\Cal P}$ into sets ${\Cal P}_0$, $\ldots$, ${\Cal P}_{J}$ 
where ${\Cal P}_0$ contains those primes in ${\Cal P}$ with $|f(p)| \le 1/2$, 
and for $1\le j\le J$ the set ${\Cal P}_j$ contains those primes 
in ${\Cal P}$ with $2^{j-2}< |f(p)| \le 2^{j-1}$.

Let $k\ge 1$ be a natural number.  Define   ${\Cal N}_0(k)$ 
to be the set of integers divisible by at most $k$ distinct {\sl squares of primes in} ${\Cal P}_0$. 
For $1\le j\le J$ we define ${\Cal N}_j(k)$ to be the set of integers 
divisible by at most $k$ distinct primes in ${\Cal P}_j$.  

\proclaim{Proposition 3.3} Keep the notations above.  For $2 \le k \le |{\Cal P}_0|/4$ we have
$$ 
\sum\Sb n\le x/y \\ n \in {\Cal N}_0(k) \endSb |f(n)|^2 \le \frac{4k}{|{\Cal P}_0|} \sum_{n\le x} |f(n)|^2.
$$ 
Further, if $1\le j\le J$ and $1\le k \le |{\Cal P}_j|/4-1$ we have
$$ 
\sum\Sb n\le x/q\\ n\in {\Cal N}_j(k)\endSb |f(n)|^2 \le \frac{2^{12} k^2}{2^{4j}|{\Cal P}_j|^2}  \sum_{n\le x} |f(n)|^2.
$$ 
\endproclaim

\demo{Proof}  Note that if $p\in {\Cal P}_0$ then 
$|f(p)|\le 1/2$, and so $|f(p^2)| = |f(p)^2 -1| \ge 3/4$.  
Therefore 
$$ 
\align
\sum\Sb n\le x/y \\ n \in {\Cal N}_0(k) \endSb |f(n)|^2 \Big( \sum\Sb p\in {\Cal P}_0 \\ p^2 \nmid n\endSb |f(p^2)|^2 
\Big) &\ge \frac {9}{16} (|{\Cal P}_0| -k ) \sum\Sb n\le x/y \\ n\in {\Cal N}_0(k) \endSb |f(n)|^2\\
&\ge \frac {27}{64} |{\Cal P}_0| \sum\Sb n\le x/y\\ n\in {\Cal N}_0(k) \endSb |f(n)|^2. \tag{3.3}\\
\endalign
$$ 
If $p \in {\Cal P}_0$ and $p^2 \nmid n$ then we claim that $|f(n) f(p^2)| \le |f(p^2n)|$.  
If $p \nmid n$ then equality holds in this claim.  If $p$ exactly divides $n$ then the 
claim amounts to $|f(p^3)| \ge |f(p)f(p^2)|$, and to see this 
note that $f(p^3)= f(p)(f(p)^2 -2)$ and $f(p^2)= f(p)^2-1$, and the 
estimate $|f(p)^2-2| \ge |f(p)^2-1|$ holds since $|f(p)|\le 1/2$.   Therefore the 
LHS of (3.3)  is 
$$ 
\le \sum\Sb m \le x \endSb |f(m)|^2\Big( \sum\Sb m= np^2\\ n\le x/y, n\in {\Cal N}_0(k) \\  p\in {\Cal P}_0, p^2 \nmid n \endSb 1 \Big) \le (k+1) \sum_{m\le x} |f(m)|^2, 
$$   
since the sum over $n$ and $p$ above is zero unless $m$ is divisible by at most $k+1$ 
squares of primes in ${\Cal P}_0$ and in this case the number of choices for $p$ in that
sum is at most $k+1$.  
Since $k\ge 2$ the stated bound follows.

The second assertion is similar.   If $p \in {\Cal P}_j$ then $|f(p)| \ge 2^{j-2}$.  Therefore 
$$ 
\align
\sum\Sb n\le x/y\\ n\in {\Cal N}_j(k)  \endSb |f(n)|^2 \Big( \sum\Sb p_1, p_2 \in {\Cal P}_j \\ p_1 <p_2 \\ p_i \nmid n\endSb 
|f(p_1p_2)|^2 \Big) &\ge 2^{4(j-2)} \binom{|{\Cal P}_j|-k}{2} \sum\Sb n\le x/y \\ n\in{\Cal N}_j(k)\endSb |f(n)|^2\\
&\ge 2^{4j-8} \frac{9|{\Cal P}_j|^2}{32} \sum\Sb n\le x/y\\ n\in {\Cal N}_j(k)\endSb |f(n)|^2.\\
\endalign
$$ 
But the LHS sums terms of the form $|f(m)|^2$ where $m=np_1p_2 \le x$ and 
$m$ is divisible by at most $k+2$ distinct primes in ${\Cal P}_j$; 
moreover each such term appears at most $\binom{k+2}{2}$ times on the LHS.  
Therefore the LHS above is 
$$ 
\le \binom{k+2}{2} \sum_{n\le x} |f(n)|^2 \le 3 k^2 \sum_{n\le x} |f(n)|^2, 
$$ 
and the Proposition follows in this case.

\enddemo 

\head 4.  Proof of Theorem 3 \endhead 

\noindent Consider the set of values $y$ with ${\Cal F}(y) \ge 10^8 \log (e y)/\sqrt{y}$.  
Pick a ``maximal" element from this set; precisely, a value $y$ belonging 
to the exceptional set, but such that no value larger than $y+1$ is in this 
set.  We shall use the work in \S 3 with this maximal value of $y$ in mind, 
and employ the notation introduced there.   The argument splits into two cases: since $|{\Cal P}| 
\ge \sqrt{y}/(2\log y)$ we must have
either $|{\Cal P}_0| \ge \sqrt{y}/(4\log y)$, or that $|{\Cal P}_j| \ge \sqrt{y}/(4J \log y)$ 
for some $1\le j\le J$. 

\subhead Case 1: $|{\Cal P}_0| \ge \sqrt{y}/(4\log y)$ \endsubhead 

Take $K=[|{\Cal P}_0| {\Cal F}(y)/8]$, so that $10^4 \le K \le |{\Cal P}_0|/4$.  Proposition 
3.3 gives that 
$$ 
\sum\Sb n\le x/y \\ n\in {\Cal N}_0 (K) \endSb |f(n)|^2 \le \frac{1}{2} {\Cal F}(y)\sum_{n\le x} |f(n)|^2,
$$ 
so that 
$$ 
\sum\Sb n\le x/y \\ n \not\in {\Cal N}_0(K) \endSb |f(n)|^2 \ge \frac 12 {\Cal F}(y) \sum_{n\le x} |f(n)|^2. 
\tag{4.1}
$$ 

If $n\not\in {\Cal N}_0(K)$ then $n$ must be divisible by at least 
$K+1$  squares of primes in ${\Cal P}_0$.   There are 
$\binom{|{\Cal P}_0|}{K+1}$ integers that are products of 
exactly $K+1$ primes from ${\Cal P}_0$.  Each of these integers exceeds $(\sqrt{y}/2)^{K+1}$, 
and a number $n\not\in {\Cal N}_0(K)$ must be divisible by the 
square of one of these integers.   
Thus, using the second bound in Proposition 3.2 we find that 
the LHS of (4.1)  is 
$$ 
\le \binom{|{\Cal P}_0|}{K+1} 3^{K+1} \cdot 3^{K+1} {\Cal F}(y(y/4)^{K+1}) \sum_{n\le x} |f(n)|^2. \tag{4.2}
$$
Since 
$$
\binom {|{\Cal P}_0|}{K+1} \le \frac{|{\Cal P}_0|^{K+1} }{(K+1)!} 
< \Big(\frac{e |{\Cal P}_0|}{K+1}\Big)^{K+1} < \Big(\frac{24}{{\Cal F}(y)}\Big)^{K+1},
$$ 
and, by the maximality of $y$, 
$$ 
{\Cal F}(y(y/4)^{K+1}) \le 10^8 \cdot 2^{K+1} \frac{1+(K+2)\log y}{y^{(K+2)/2}} < 10^{8} \cdot 2^{K+1} 
\cdot (10^{-8} {\Cal F}(y))^{K+2},
$$ 
we deduce that 
the quantity in (4.2) is 
$$ 
\le \Big( \frac{432}{10^{8}}\Big)^{K+1} {\Cal F}(y) \sum_{n\le x} |f(n)|^2 < 
\frac 12 {\Cal F}(y) \sum_{n\le x} |f(n)|^2, 
$$
which contradicts (4.1).  This completes our argument for the first case. 

 \subhead Case 2:  $|{\Cal P}_j | \ge \sqrt{y}/(4J \log y)$ for some $1\le j\le J$ \endsubhead

Here we take $K= [2^{2j-9} |{\Cal P}_j| {\Cal F}(y)^{\frac 12}]$.   Using that 
$J \le 3 +(\log (1/{\Cal F}(y)))/(4\log 2) \le (\log y)/4$, and $y\ge 10^{16}$ we may check 
that $K \ge 10$.  Moreover, for any prime $p$ in ${\Cal P}_j$ we have 
$2^{2j-4} \le |f(p)|^2 \le 4/{\Cal F}(y)^{\frac 12}$ by Lemma 3.1, and so 
$K\le |{\Cal P}_j|/8$.  Thus the second part of Proposition 3.3 applies, and 
it shows that 
$$ 
\sum\Sb n\le x/y\\ n \in {\Cal N}_j(K)\endSb |f(n)|^2 \le \frac 12 {\Cal F}(y) \sum_{n\le x} |f(n)|^2. 
$$
Therefore
$$
\sum\Sb n\le x/y \\ n\not\in {\Cal N}_j(K) \endSb |f(n)|^2 \ge \frac 12 {\Cal F}(y) \sum_{n\le x} |f(n)|^2. 
\tag{4.3}
$$

If $n\not\in {\Cal N}_j(K)$ then $n$ must be divisible by one of the $\binom{|{\Cal P}_j|}{K+1}$ 
integers composed of exactly $K+1$ primes in ${\Cal P_j}$.  Each of those numbers 
exceeds $(\sqrt{y}/2)^{K+1}$.  Appealing to the first part of Proposition 3.2 we 
find that the LHS of (4.3)  is 
$$ 
 \binom{|{\Cal P}_j|}{ K+1} 2^{2j(K+1)} {\Cal F}(y(\sqrt{y}/2)^{K+1}) \sum_{n\le x } |f(n)|^2. \tag{4.4}
$$
 Since
 $$ 
 \binom{|{\Cal P}_j|}{K+1} \le \frac{|{\Cal P}_j|^{K+1}}{(K+1)!} \le \Big( 
 \frac{e |{\Cal P}_j|}{K+1}\Big)^{K+1} < \Big( \frac{2^{11}}{2^{2j} {\Cal F}(y)^{\frac 12}}\Big)^{K+1}, 
 $$ 
 and, by the maximality of $y$, 
 $$ 
 {\Cal F}(y (\sqrt{y}/2)^{K+1} ) \le 10^8 \cdot 2^{\frac{K+1}{2}} \frac{1+\frac{K+3}{2}\log y}{y^{\frac {K+3}{2}}} 
 < 10^8 \cdot 2^{\frac{K+1}{2}} ( 10^{-8} {\Cal F}(y))^{\frac{K+3}{2}},
 $$ 
 we deduce that the quantity in (4.4) is 
 $$ 
 \le \Big( \frac{2^{23}}{10^8} \Big)^{\frac{K+1}{2}} {\Cal F}(y) \sum_{n\le x} |f(n)|^2 
 < \frac 12 {\Cal F}(y) \sum_{n\le x} |f(n)|^2,
 $$ 
 which contradicts (4.3).  This completes our argument in the second case, and 
 hence also the proof of Theorem 3.

 \enddocument